\documentclass[12pt]{article}
\usepackage{amssymb,amsmath}
\usepackage{mathrsfs}
\usepackage{graphicx}

\newtheorem{theorem}{Theorem}
\newtheorem{lemma}{Lemma}

\begin{document}

\begin{center}
\large{\bf ON THE MAXIMAL DISPLACEMENT\\
OF CATALYTIC BRANCHING RANDOM WALK}
\end{center}
\vskip0,5cm
\begin{center}
Ekaterina Vl. Bulinskaya\footnote{ \emph{Email address:} {\tt
bulinskaya@yandex.ru}}$^,$\footnote{The work is supported by Russian Science Foundation under grant 17-11-01173-Ext and is fulfilled
at Novosibirsk State University. The author is Associate Professor of the Lomonosov Moscow State University.}
\end{center}
\vskip1cm

\begin{abstract}
We study the distribution of the maximal displacement of particles positions for the whole time of the population existence in the model of critical and subcritical catalytic branching random walk on $\mathbb{Z}$. In particular, we prove that in the case of simple symmetric random walk on $\mathbb{Z}$, the distribution of the maximal displacement has ``a heavy tail'' decreasing as a function of the power $1/2$ or $1$, when the branching process is critical or subcritical, respectively. These statements describe new effects which do not arise in the corresponding investigations of the maximal displacement of critical and subcritical branching random walks on $\mathbb{Z}$.


\vskip0,5cm 2010 {\it AMS classification}: 60J80, 60F05.

\end{abstract}

\section{Introduction}

Problems of the rate of population propagation (of particles, bacteria, individuals, genes etc.) in space have been attracting attention of researchers for a long time. It suffices to indicate, for example, the survey \cite{Shi_LNM_15} and the paper \cite{Ray_Hazra_Roy_Soulier_19} devoted to branching random walk (BRW), recent works \cite{Oz_19} and \cite{Shiozawa_19} in which branching Brownian motion (BBM) is investigated. Among the models describing the evolution of popu\-lation in space a special place is occupied by catalytic branching processes and, in particular, catalytic branching Brownian motion (CBBM), see, e.g., publica\-tions \cite{Bocharov_Wang_19} and \cite{Wang_Zong_19}. We note also a catalytic branching random walk (CBRW), the present work focusing attention on its study. A distinctive feature of catalytic branching processes is furnishing the space with catalysts and only there a particle may produce offspring or die. Outside the catalysts particles can only move in space. Thus, a particle evolution depends on its spatial location.

Till now the problems of propagation of particle population have been considered in the case of {\it supercritical} CBRW on $\mathbb{Z}^d$, $d\in\mathbb{N}$, see, e.g., papers \cite{Molchanov_Yarovaya_12}, \cite{Carmona_Hu_14} and \cite{B_UMN_19}. In the supercritical regime the particle population in CBRW survives with positive probability, and in the case of survival the total and local numbers of particles grow exponentially-fast in course of time (see \cite{B_TPA_15} and \cite{B_Doklady_15}). While in {\it critical} and {\it subcritical} regimes the population degenerates locally with probability 1, in some cases it can survive globally with positive probability (for global and local extinction see, e.g., \cite{Bertacchi_Zucca_14}). Therefore, in supercritical CBRW the rate of population propagation, as time grows unboundedly, is of interest, whereas in critical and subcritical regimes the main attention is paid to the maximal displacement of particles for the whole history of the population existence.

It turned out that the rate of propagation of particle population in CBRW depends essentially on ``heaviness'' of the distribution tails of the walk jump. For this reason in the series of papers \cite{B_SPA_18}, \cite{B_Arxiv_18} and \cite{B_Arxiv_19} we had to consider separately the cases of ``light tails'', regularly varying tails and that of semi-exponential distribution of the walk jump. In the present work we are interested in critical and subcritical CBRW on $\mathbb{Z}$. Thus, in the context of investigation of population propagation, the aim of the work consists in the study of the maximal displacement of particles for the whole history of population existence.

For the distribution functions of the maximal displacement of particles for the whole history of the process we derive a system of equations having a unique solution. In the system there arise the probabilities related to the behavior of the random walk only on the time-interval from the exit moment from a catalyst till the moment of the first return to it or the moment of the first hitting another catalyst. Such probabilities have not been studied previously for arbitrary random walks. However, in particular but important case of a simple random walk (i.e. the jumps of the walk are performed to the nearest-neighbor points of the lattice $\mathbb{Z}$) these probabilities can be found on the basis of solution to the classical ``ruin'' problem. Therefore, in this case we investigate asymptotic behavior of the distribution tails of the maximal displacement of particles in critical and subcritical CBRW on $\mathbb{Z}$. Whenever the simple random walk has a drift, the obtained results are natural and do not surprise: the distribution tail of the maximal displacement either decays exponentially-fast or the random variable under conside\-ration is an extended one. However, whenever the simple random walk is symmetric, the new results appear unexpected and radically differ from the known statements for BRW studied in the papers \cite{Lalley_Shao_15} and \cite{Neuman_Zheng_17}.

The studies of CBRW different from estimation of the rate of population propagation can be found, e.g., in works \cite{Khristolyubov_Yarovaya_19}, \cite{Platonova_Ryadovkin_17} and \cite{VT_13}.

\section{Main results}

The description of CBRW model on $\mathbb{Z}$ with $N$ catalysts, forming the set $W=\{w_1,\ldots,w_N\}\subset\mathbb{Z}$, is present, e.g., in papers \cite{Carmona_Hu_14}, \cite{B_TPA_15}, \cite{B_Arxiv_18}. However, for the sake of the reader's convenience we recall it here. Assume that all the random variables are specified on a complete probability space $(\Omega,\mathcal{F},{\sf P})$, where $\Omega$ is a sample space consisting of outcomes $\omega$. Moreover, index $x$ of the probability ${\sf P}_x$ and mathematical expectation ${\sf E}_x$ marks the starting point of CBRW or of the random walk, depending on the context.

Let at the time moment $t=0$ there be a single particle on the lattice located at point $z\in\mathbb{Z}$. Whenever $z\notin W$, the movement of the particle until the first hitting the catalysts set $W$ is determined by the Markov chain $S=\{S(t),t\geq0\}$. The space-homogeneous random walk $S$ is specified by the infinitesimal matrix $Q=\left(q(x,y)\right)_{x,y\in\mathbb{Z}}$, which is assumed irreducible and conservative, i.e.
\begin{equation}\label{condition1}
q(x,y)=q(x-y,0)=q(0,y-x)\quad\mbox{and}\quad\sum\limits_{y\in\mathbb{Z}}{q(x,y)}=0,
\end{equation}
where $q(x,y)\geq0$, for $x\neq y$, and $q(x,x)\in(-\infty,0)$, for all $x,y\in\mathbb{Z}$. Whenever $z\in W$ or the particle hits the set $W$ for the first time, for instance, at a catalyst $w_k$, $k=1,\ldots,N$, then the particle spends there random time distributed exponentially with parameter $\beta_k>0$. Afterwards with probability $\alpha_k\in[0,1)$ it instantly produces a random number $\xi_k$ of offsprings located also at $w_k$ and dies. Otherwise the particle performs a jump at point $y$ with probability $-(1-\alpha_k)q(w_k,y)/q(w_k,w_k)$, $y\in\mathbb{Z}$, $y\neq w_k$, and continues its walking until the next hitting the catalysts set. The new particles behave as independent copies of the parent particle.

Denote by $f_k(s):={\sf E}{s^{\xi_k}}$, $s\in[0,1]$, the probability generating function of the random variable $\xi_k$ and set $m_k:={\sf E}{\xi_k}=f_k'(1)<\infty$, $k=1,\ldots,N$. We exclude the deterministic case when $f_k(s)=s$, $s\in[0,1]$, for all $k=1,\ldots,N$.

In paper~\cite{B_TPA_15} there was proposed a classification such that CBRW is called super\-critical, critical or subcritical, whenever the Perron root (i.e. the maximal positive eigenvalue) of matrix
\begin{equation}\label{criticality_CBRW}
D=\left(\delta_{i,j}\alpha_im_i+(1-\alpha_i){_{W_j}}F_{w_i,w_j}(\infty)\right)_{i,j=1}^N
\end{equation}
is larger than $1$, equal to $1$ or less than $1$, respectively. Here $\delta_{i,j}=1$, if $i=j$, and $\delta_{i,j}=0$ otherwise. Let also $W_j:=W\setminus\left\{w_j\right\}$. Then ${_{W_j}}F_{w_i,w_j}(\infty)$ is a probability of hitting point $w_j$ by the random walk avoiding the set $W_j$, whenever the starting point is $w_i$, $i,j=1,\ldots,N$. In paper~\cite{B_Doklady_15} there was established that in supercritical CBRW only the total and local numbers of particles grow exponentially-fast over time, whereas the probabilities of global and local survival are positive. The rate of propagation of particle population in supercritical CBRW was studied in works~\cite{Molchanov_Yarovaya_12}, \cite{Carmona_Hu_14}, \cite{B_SPA_18}--\cite{B_UMN_19}. Since, contrary to supercritical CBRW, in critical and subcritical CBRW the population of particles degenerates locally, it makes no sense to talk about the rate of population propagation. However one can pose a question how remote points are visited by the particles \emph{during the whole history of the population existence}. Our work is devoted to answering this question in cases of critical and subcritical CBRW on $\mathbb{Z}$.

Let $Z(t)$ be a random set of particles existing in CBRW at time moment $t\geq0$. For a particle $v\in Z(t)$, denote by $X^v(t)$ its location at time $t\geq0$. Let $M_t:=\max\{X^v(t),v\in Z(t)\}$ be the maximum of CBRW at time $t\geq0$, i.e. the location of the right-most particle existing in CBRW at time $t$. We will be interested in the random variable $M:=\max\{M_t,t\geq0\}$ being the maximal displacement (to the right from the origin) of CBRW for the whole history of the particle population. Clearly, $M\geq z$.

In formulations of Theorems~\ref{T:critical_simple_symmetric}--\ref{T:subcritical_simple_asymmetric} we consider a simple random walk $S$ on the lattice $\mathbb{Z}$. It means that
$$
\frac{q(x,x+1)}{-q(x,x)}=p,\quad \frac{q(x,x-1)}{-q(x,x)}=q,\quad q(x,y)=0,\quad\mbox{for}\quad|x-y|\geq 2,
$$
where $p+q=1$ and $p,q\in(0,1)$. Such a random walk is called symmetric, whenever $p=q$, and asymmetric otherwise. In other words, in one jump a particle performing a simple random walk on $\mathbb{Z}$ moves to the nearest point to the right with probability $p$ and to the nearest point to the left with probability $q$. A simple random walk on $\mathbb{Z}$ is recurrent if and only if it is symmetric (see, e.g., \cite{Borovkov_book_13}, Theorem~13.3.1).

To prove Theorems~\ref{T:critical_simple_symmetric}--\ref{T:subcritical_simple_asymmetric} we derive equations~(\ref{P0(M>x)_main_system})--(\ref{P0(M>x)_main_equation_arbitrary_z}) for the probabilities under consideration. These equations are valid for an arbitrary number of catalysts and for any random walk satisfying condition~(\ref{condition1}) (not only for a simple random walk). However, for the subsequent study of the solutions to the equations we have to know such properties of the random walks which can be established easily in the case of the simple random walk and constitute a separate research work in another case. Therefore, in the present work our main results are based on the assumption of the random walk simplicity.

In Theorems~\ref{T:critical_simple_symmetric}--\ref{T:subcritical_simple_asymmetric} we also assume that the set $W$ consists of a single catalyst located at the origin $0$, whereas the starting point is also located at $0$. The asymptotic results in Theorems~\ref{T:critical_simple_symmetric}--\ref{T:subcritical_simple_asymmetric} hold true under wider assumptions of any finite number of catalysts and an arbitrary starting point. The difference consists in the constants arising in the asymptotics. However, the form of these constants depends essentially on the relative location of both the starting point and the catalysts as well as on the distances between them. That is why the corresponding bulky results are not reproduced here.

In the next theorem we establish the asymptotic behavior of the distribution tail of the random variable $M$ for a \emph{critical} CBRW on $\mathbb{Z}$, in which the random walk is a simple and symmetric one. Here and after, whenever we talk about a single catalyst, we assume that, without loss of generality, it is located at $0$, and the index $1$ of the symbols $\alpha_1$, $\xi_1$, $f_1$ and $m_1$ is omitted. Since, as noted above, a simple symmetric random walk is recurrent, the probability of return from $0$ to $0$ denoted above as ${_\varnothing}F_{0,0}(\infty)$ equals $1$. Thus, the definition of a critical CBRW (see formula (\ref{criticality_CBRW})) results in the equality $\alpha m+(1-\alpha){_\varnothing}F_{0,0}(\infty)=1$ which is equivalent to $m=1$. In other words, for a recurrent random walk, CBRW with a single catalyst is critical if and only if the Galton-Watson branching process with an offspring number $\xi_1$ is critical.

\begin{theorem}\label{T:critical_simple_symmetric}
Let $f'(1)=1$ and $f''(1)=\sigma^2\in(0,\infty)$, for CBRW on $\mathbb{Z}$, in which the random walk $S$ is simple and symmetric. Then
\begin{equation}\label{P0(M>x)_asymptotics_critical_simple_symmetric}
{\sf P}_0\left(M>x\right)\sim\frac{\sqrt{1-\alpha}}{\sqrt{\alpha\sigma^2}\sqrt{x}},\quad x\to\infty.
\end{equation}
\end{theorem}

The result of Theorem~\ref{T:critical_simple_symmetric} is a counterpart of the main result of the paper \cite{Lalley_Shao_15} derived for the model of a critical BRW on $\mathbb{Z}$. However in the latter model the decay rate of the probability ${\sf P}_0\left(M>x\right)$ has an order $1/x^2$, as $x\to\infty$. Therefore, the particles in the critical CBRW manage to go father away from the origin before returning to it and, possibly, dying, than in the model of BRW, in which the particles may die at any point.

Theorem~\ref{T:subcritical_simple_symmetric} gives the solution to the same problem as in Theorem~\ref{T:critical_simple_symmetric}. The only difference is that now we consider a \emph{subcritical} CBRW on $\mathbb{Z}$.

\begin{theorem}\label{T:subcritical_simple_symmetric}
Let $m=f'(1)<1$ for a CBRW on $\mathbb{Z}$, in which the random walk $S$ is a simple and symmetric one. Then
\begin{equation}\label{P0(M>x)_asymptotics_subcritical_simple_symmetric}
{\sf P}_0\left(M>x\right)\sim\frac{1-\alpha}{2\alpha(1-m)x},\quad x\to\infty.
\end{equation}
\end{theorem}

The result of Theorem~\ref{T:subcritical_simple_symmetric} is a counterpart of the main result of the paper~\cite{Neuman_Zheng_17} devoted to a subcritical BRW on $\mathbb{Z}$. However, in the latter case the probability ${\sf P}_0\left(M>x\right)$ decays exponentially-fast that differs importantly from our result. This difference is connected again with possible dying of the particles at any point of the lattice in the BRW model.

Theorems~\ref{T:critical_simple_symmetric} and \ref{T:subcritical_simple_symmetric} are focused on the case of a simple \emph{symmetric} random walk on $\mathbb{Z}$. The two following theorems are devoted to investigation of critical and subcritical CBRW in which the random walk is a simple and \emph{asymmetric} one, i.e. it has a drift to the right, whenever $p>q$, or to the left, whenever $p<q$. Because of a drift the random walk is no longer a recurrent one. Correspondingly, the criticality condition of CBRW changes as well. Now ${r:=1-{_\varnothing}F_{0,0}(\infty)\in(0,1)}$, and according to (\ref{criticality_CBRW}) the criticality of CBRW implies that $\alpha m+(1-\alpha)(1-r)=1$ which is equivalent to $m=1+r\alpha^{-1}(1-\alpha)$.

In the next theorem we estimate the distribution tail of the random variable $M$ for a \emph{critical} CBRW on $\mathbb{Z}$, in which the underlying random walk is \emph{simple and asymmetric}.

\begin{theorem}\label{T:critical_simple_asymmetric}
Let $m=1+r\alpha^{-1}(1-\alpha)$ and $f''(1)=\sigma^2\in(0,\infty)$ for a CBRW on $\mathbb{Z}$, in which the random walk is a simple and asymmetric one. Then the following relations hold
\begin{equation}\label{P0(M>x)_asymptotics_critical_simple_asymmetric_p<q}
{\sf P}_0\left(M>x\right)\sim\frac{\sqrt{2(1-\alpha)(q-p)}}{\sqrt{\alpha\sigma^2}}\left(\frac{p}{q}\right)^{\frac{x+1}{2}},\quad\mbox{whenever}\quad p<q,
\end{equation}
\begin{equation}\label{P0(M>x)_asymptotics_critical_simple_asymmetric_p>q}
{\sf P}_0\left(M>x\right)\to s_0,\quad\mbox{whenever}\quad p>q,
\end{equation}
as $x\to\infty$, where $s_0\in(0,1)$ is a unique solution to equation $\alpha(1-f(1-s))+\left(2q(1-\alpha)-1\right)s+(1-\alpha)(p-q)=0$ with respect to unknown variable $s$, $s\in[0,1]$.
\end{theorem}

The following result contains solution to the same problem which is the subject of Theorem~\ref{T:critical_simple_asymmetric}, but now for \emph{subcritical} CBRW on $\mathbb{Z}$.

\begin{theorem}\label{T:subcritical_simple_asymmetric}
Let $m<1+r\alpha^{-1}(1-\alpha)$, for CBRW on $\mathbb{Z}$, in which the random walk is simple and asymmetric. Then
\begin{equation}\label{P0(M>x)_asymptotics_subcritical_simple_asymmetric_p<q}
{\sf P}_0\left(M>x\right)\sim\frac{(1-\alpha)(q-p)}{1-2p(1-\alpha)-\alpha m}\left(\frac{p}{q}\right)^{x+1},\quad\mbox{whenever}\quad p<q,
\end{equation}
\begin{equation}\label{P0(M>x)_asymptotics_subcritical_simple_asymmetric_p>q}
{\sf P}_0\left(M>x\right)\to s_0,\quad\mbox{whenever}\quad p>q,
\end{equation}
as $x\to\infty$, where $s_0\in(0,1)$ is a unique root of equation $\alpha(1-f(1-s))+\left(2q(1-\alpha)-1\right)s+(1-\alpha)(p-q)=0$ with respect to unknown variable $s$, $s\in[0,1]$.
\end{theorem}

The results of Theorems~\ref{T:critical_simple_asymmetric} and \ref{T:subcritical_simple_asymmetric} are expected. Namely, if the random walk $S$ has a drift to the left ($p<q$), then the particles in CBRW do not manage to go far away to the right, since they drift to the left. Conversely, if the random walk $S$ has a drift to the right ($p>q$), then there are particles in CBRW which will go away to the right to ``infinity'', and therefore $M=\infty$ with positive probability $s_0$.

Thus, in the case of a simple random walk we find the asymptotic behavior of probability ${\sf P}_0\left(M>x\right)$, as $x\to\infty$, in critical and subcritical CBRW on $\mathbb{Z}$ with a single catalyst at $0$. The formulated Theorems~\ref{T:critical_simple_asymmetric} and \ref{T:subcritical_simple_asymmetric} in the case of asymmetric simple random walk are not surprising and are presented for completeness of the picture. The results of Theorems~\ref{T:critical_simple_symmetric} and \ref{T:subcritical_simple_symmetric} describe new effects and are of the main interest. Indeed, they are radically different from the corresponding statements for BRW on $\mathbb{Z}$ studied in \cite{Lalley_Shao_15} and \cite{Neuman_Zheng_17}. The results obtained by us are the first investigation in the domain of description of population propagation in critical and subcritical CBRW. It is worthwhile to note that visible differences in the propagation of particle population in supercritical CBRW and supercritical BRW were revealed only in the second term of the asymptotic expansions for their corresponding maximums (see, e.g., \cite{B_TPA_19}, \cite{Carmona_Hu_14} and \cite{Mallein_16}). Meanwhile, as shown in our investigations, in critical and subcritical CBRWs and the corresponding critical and subcritical BRWs the differences are noticeable already in the first asymptotic approximation of the probability ${\sf P}_0\left(M>x\right)$, as $x\to\infty$.

\section{Proofs}

First of all, recall the definition (see, e.g., \cite{B_SPL_14}) of hitting times under taboo which we need for deriving equations with respect to the probability ${\sf P}_z(M>x)$, $z\in\mathbb{Z}$, under consideration. Set
$$
\tau_x:=\inf\{t>0:S(t)\neq S(0)\}\mathbb{I}\left(S(0)=x\right),
$$
i.e. introduce the exit moment of the random walk $S$ from the starting point $x\in\mathbb{Z}$. As usual, $\mathbb{I}\left(A\right)$ is an indicator of event $A\in\mathcal{F}$. Denote by
$$
{_H}\tau_{x,y}:=\inf\left\{t\geq\tau_x:S(t)=y,S(u)\notin H,\tau\leq u\leq t\right\}\mathbb{I}\left(S(0)=x\right)
$$
time of (the first) hitting by the random walk $S$ the point $y\in\mathbb{Z}$ under taboo on the visit of set $H\subset\mathbb{Z}$, $y\notin H$, when the walk starts at point $x\in\mathbb{Z}$. Whenever the trajectory of the random walk $S(\cdot,\omega)$ after the start at point $x$ visits the set $H$ before hitting point $y$, then we naturally set ${_H}\tau_{x,y}(\omega)=\infty$. Note that ${_{W_j}}F_{w_i,w_j}(\infty)={\sf P}_{w_i}\left({_{W_j}}\tau_{w_i,w_j}<\infty\right)$.

\begin{lemma}\label{L:main_equations}
The following system of equations holds true with respect to probabilities ${\sf P}_{w_i}\left(M>x\right)$, $x\in\mathbb{Z}$, $i=1,\ldots,N$:
\begin{eqnarray}\label{P0(M>x)_main_system}
& &\!\!{\sf P}_{w_i}\left(M>x\right)=\alpha_i\left(1-f_i\left(1-{\sf P}_{w_i}\left(M>x\right)\right)\right)\\
&+&\!\!(1-\alpha_i)\!\sum_{j=1}^N\!{\sf P}_{w_i}\!\left(\max\left\{S(t),0\leq t\leq{_{W_j}\tau_{w_i,w_j}}\right\}\leq x,\,{_{W_j}\tau_{w_i,w_j}}<\infty\right)\!{\sf P}_{w_j}\!\left(M\!>\!x\right)\nonumber\\
&+&\!\!(1-\alpha_i)\,{\sf P}_{w_i}\!\left(\max\left\{S(t),0\leq t\leq\min\limits_{j=1,\ldots,N}{_{W_j}\tau_{w_i,w_j}}\right\}>x\right),\nonumber
\end{eqnarray}
where $W_j:=W\setminus\{w_j\}$, $j=1,\ldots,N$.

The case of the start CBRW at an arbitrary point $z\in\mathbb{Z}\setminus W$ is reduced to the previous one:
\begin{eqnarray}\label{P0(M>x)_main_system_arbitrary_z}
& &{\sf P}_z\left(M>x\right)={\sf P}_z\left(\max\left\{S(t),0\leq t\leq\min\limits_{i=1,\ldots,N}{_{W_i}\tau_{z,w_i}}\right\}>x\right)\\
\!\!&+&\!\!\sum_{i=1}^N{\sf P}_z\left(\max\left\{S(t),0\leq t\leq{_{W_i}\tau_{z,w_i}}\right\}\leq x,\,{_{W_i}\tau_{z,w_i}}<\infty\right){\sf P}_{w_i}\!\left(M>x\right),\nonumber
\end{eqnarray}
where, evidently, ${\sf P}_z\left(M>x\right)=1$, for $x<z$.

In particular, whenever $W=\{0\}$, the system of equations (\ref{P0(M>x)_main_system}) transforms into the following equation with respect to ${\sf P}_0\left(M>x\right)$:
\begin{eqnarray}\label{P0(M>x)_main_equation}
& &{\sf P}_0\left(M>x\right)=\alpha\left(1-f\left(1-{\sf P}_0\left(M>x\right)\right)\right)\\
&+&(1-\alpha){\sf P}_0\left(\max\left\{S(t), 0\leq t\leq\tau_{0,0}\right\}\leq x,\tau_{0,0}<\infty\right){\sf P}_0\left(M>x\right)\nonumber\\
&+&(1-\alpha){\sf P}_0\left(\max\left\{S(t),0\leq t\leq\tau_{0,0}\right\}>x\right).\nonumber
\end{eqnarray}
The case of the start at point $z\neq0$, $z\in\mathbb{Z}$, is reduced to the previous one as well:
\begin{eqnarray}\label{P0(M>x)_main_equation_arbitrary_z}
{\sf P}_z\left(M>x\right)&=&{\sf P}_z\left(\max\left\{S(t),0\leq t\leq\tau_{z,0}\right\}>x\right)\\
&+&{\sf P}_z\left(\max\left\{S(t),0\leq t\leq\tau_{z,0}\right\}\leq x,\tau_{z,0}<\infty\right){\sf P}_0\left(M>x\right).\nonumber
\end{eqnarray}
The system (\ref{P0(M>x)_main_system}) and, in particular, equation (\ref{P0(M>x)_main_equation}) have a unique solution on the intervals $[0,1]^N$ and $[0,1]$, respectively.
\end{lemma}
{\sc Proof.} To reduce the volume of the work we consider the most illustrative case $W=\{0\}$ and $z=0$. The rest of the proof of Lemma~\ref{L:main_equations} is conducted on the basis of the same ideas as this main case. In view of the formula of total probability and according to the description of the CBRW model we have
\begin{eqnarray*}
{\sf P}_0(M\leq x)&=&\alpha\sum_{k=0}^{\infty}{\sf P}(\xi=k)\left({\sf P}_0(M\leq x)\right)^k\\
&+&(1-\alpha){\sf P}_0\left(\max\left\{S(t),0\leq t\leq\tau_{0,0}\right\}\leq x,\tau_{0,0}<\infty\right){\sf P}_0\left(M\leq x\right)\\
&+&(1-\alpha){\sf P}_0\left(\tau_{0,0}=\infty, S(t)\leq x,t\geq0\right),
\end{eqnarray*}
which is equivalent to (\ref{P0(M>x)_main_equation}).

The solution to equation (\ref{P0(M>x)_main_equation}) with respect to ${\sf P}_0\left(M>x\right)$ always exists and is unique, since the solution to equation
$$
\alpha(1-f(1-s))=s\left(1-(1-\alpha)p_1\right)-(1-\alpha)p_2
$$
exists and is unique, for $s\in[0,1]$, where
\begin{eqnarray*}
p_1&:=&{\sf P}_0\left(\max\left\{S(t),0\leq t\leq\tau_{0,0}\right\}\leq x,\tau_{0,0}<\infty\right),\\
p_2&:=&{\sf P}_0\left(\max\left\{S(t),0\leq t\leq\tau_{0,0}\right\}>x\right),
\end{eqnarray*}
and, obviously, $p_1+p_2\leq1$. Indeed , we have $0=\alpha(1-f(1))>-(1-\alpha)p_2$, for $s=0$, and $\alpha(1-f(0))\leq\alpha\leq1-(1-\alpha)p_1-(1-\alpha)p_2$, for $s=1$. Therefore, whenever at least one inequality in the latter relation is strict, then the graphs of functions $\alpha(1-f(1-s))$ and $s\left(1-(1-\alpha)p_1\right)-(1-\alpha)p_2$, for $s\in[0,1]$, have a unique (by the convexity of function $f$) intersection point on the interval $(0,1)$. Whenever $\alpha(1-f(0))=\alpha=1-(1-\alpha)p_1-(1-\alpha)p_2$ (it is possible in the case of a recurrent random walk and the null probability of a particle death without giving offspring), then the intersection of the mentioned graphs is at point $s=1$ and there are no other intersection points, since $\frac{d}{ds}\left(\alpha(1-f(1-s))\right)_{s=1}<\alpha<\frac{d}{ds}\left(s\left(1-(1-\alpha)p_1\right)-(1-\alpha)p_2\right)_{s=1}$. $\square$

It follows from (\ref{P0(M>x)_main_equation}) that the asymptotic behavior of ${\sf P}_0\left(M>x\right)$, as $x\to\infty$, is determined by that of ${\sf P}_0\left(\max\left\{S(t),0\leq t\leq\tau_{0,0}\right\}\leq x,\tau_{0,0}<\infty\right)$ and ${\sf P}_0\left(\max\left\{S(t),0\leq t\leq\tau_{0,0}\right\}>x\right)$. Under general assumptions on the random walk these probabilities have not been studied. However, in a particular but important case of the simple random walk the investigation of these probabilities can be reduced to the already solved classic ``ruin problem''. In the following two lemmas the formulae for these probabilities are derived separately for the cases of a simple symmetric and a simple asymmetric random walk.

\begin{lemma}\label{L:P_0(max>x)_simple_symmetric}
For a simple symmetric random walk $S$ on $\mathbb{Z}$ and $x\in\mathbb{N}$, the following equalities hold true
\begin{equation}\label{P_0(S(t)<=x)_simple_symmetric}
{\sf P}_0\left(\max\left\{S(t),0\leq t\leq\tau_{0,0}\right\}\leq x,\tau_{0,0}<\infty\right)=\frac{2x+1}{2(x+1)},
\end{equation}
\begin{equation}\label{P_0(max>x)_simple_symmetric}
{\sf P}_0\left(\max\left\{S(t),0\leq t\leq\tau_{0,0}\right\}>x\right)=\frac{1}{2(x+1)}.
\end{equation}
\end{lemma}
{\sc Proof.}
A simple symmetric random walk on $\mathbb{Z}$ is a recurrent one (see, e.g., \cite{Borovkov_book_13}, Theorem~13.3.1). Consequently, $\tau_{0,0}=\infty$ with probability $0$ and
$$
{\sf P}_0\left(\max\left\{S(t),0\leq t\leq\tau_{0,0}\right\}\leq x,\tau_{0,0}<\infty\right)=1-{\sf P}_0\left(\max\left\{S(t),0\leq t\leq\tau_{0,0}\right\}>x\right).
$$

Let us derive a formula for ${\sf P}_0\left(\max\left\{S(t),0\leq t\leq\tau_{0,0}\right\}>x\right)$ in the case of a simple symmetric random walk. Since the jumps of the random walk may occur to the neighbor points, then in the random event $\{\omega:\max\left\{S(t),0\leq t\leq\tau_{0,0}\right\}>x\}$ there are only the trajectories of $S$ which from starting point $0$ pass to point $1$ and then hit point $x+1$ before point $0$. Thus, taking into account the results of the classic ``ruin problem'' (see, e.g., \cite{Shiryaev_Probability1_16}, Ch.~1, \S9, formula~(14)) we come to relation (\ref{P_0(max>x)_simple_symmetric}) and, therefore, to relation (\ref{P_0(S(t)<=x)_simple_symmetric}). Lemma~\ref{L:P_0(max>x)_simple_symmetric} is proved completely. $\square$

\begin{lemma}\label{L:P_0(max>x)_simple_asymmetric}
For a simple asymmetric random walk on $\mathbb{Z}$, the following formulae are valid:
\begin{equation}\label{P_0(S(t)<=x)_simple_asymmetric}
{\sf P}_0\left(\max\left\{S(t),0\leq t\leq\tau_{0,0}\right\}\leq x,\tau_{0,0}<\infty\right)=p\frac{\left(q/p\right)^{x+1}-(q/p)}{\left(q/p\right)^{x+1}-1}+\min\{p,q\},
\end{equation}
\begin{equation}\label{P_0(max>x)_simple_asymmetric}
{\sf P}_0\left(\max\left\{S(t),0\leq t\leq\tau_{0,0}\right\}>x\right)=\frac{q-p}{\left(q/p\right)^{x+1}-1},
\end{equation}
for each $x\in\mathbb{N}$.
\end{lemma}
{\sc Proof.} According to the total probability formula we have
\begin{eqnarray}\label{P_0(max>x)_simple_asymmetric_=}
&&{\sf P}_0\!\left(\max\left\{S(t),0\!\leq\! t\!\leq\!\tau_{0,0}\right\}\!\leq\! x,\tau_{0,0}\!<\!\infty\right)\!=\!{\sf P}_0\!\left(S(t)\!<\!0,\tau_0\!\leq\! t\!<\!\tau_{0,0},\tau_{0,0}\!<\!\infty\right)\\
&+&{\sf P}_0\left(S(t)\in(0,x],\tau_0\leq t<\tau_{0,0},\tau_{0,0}<\infty\right)\nonumber\\
&=&q\,{\sf P}_{-1}\left(\exists t_2:S(t_2)=0,S(t)\neq0,0\leq t<t_2\right)\nonumber\\
&+&p\,{\sf P}_1\left(\exists t_1:S(t_1)=0,S(t)\neq0,S(t)\neq x+1,0\leq t<t_1\right).\nonumber
\end{eqnarray}
Here ${\sf P}_1\left(\exists t_1:S(t_1)=0,S(t)\neq0,S(t)\neq x+1,0\leq t<t_1\right)$ is a probability of an exit of the random walk $S$ from a stripe $(0,x+1)$ through the lower boundary, when the starting point is located at $1$, see \cite{Shiryaev_Probability1_16}, Ch.~1, \S9, formula~(13). Similarly ${\sf P}_{-1}\left(\exists t_2:S(t_2)=0,S(t)\neq0,0\leq t<t_2\right)$ is a probability of an exit of the random walk $S$ from a stripe $(-\infty,0)$ through the upper boundary, when the starting point is located at $-1$. The latter probability can be found due to formula~(13) in \cite{Shiryaev_Probability1_16}, Ch.~1, \S9, as well, but now the probabilities $p$ and $q$ should be swapped and let the upper boundary $B$ tend to infinity. Returning to representation~(\ref{P_0(max>x)_simple_asymmetric_=}) and substituting the found expressions for the probabilities, we get
$$
{\sf P}_0\left(\max\left\{S(t),0\leq t\leq\tau_{0,0}\right\}\leq x,\tau_{0,0}<\infty\right)=p\frac{(q/p)^{x+1}-(q/p)}{(q/p)^{x+1}-1}+q\min\left\{\frac{p}{q},1\right\},
$$
which coincides with relation (\ref{P_0(S(t)<=x)_simple_asymmetric}).

In a similar way, with the help of formula~(10) in \cite{Shiryaev_Probability1_16}, Ch.~1, \S9, we obtain
\begin{eqnarray*}
& &{\sf P}_0\left(\max\left\{S(t),0\leq t\leq\tau_{0,0}\right\}>x\right)\\
&=&p\,{\sf P}_1\left(\exists t_3:S(t_3)=x+1,S(t)\neq0,S(t)\neq x+1,0\leq t<t_3\right)=p\frac{(q/p)-1}{(q/p)^{x+1}-1}.
\end{eqnarray*}
Lemma~\ref{L:P_0(max>x)_simple_asymmetric} is proved completely. $\square$

Let us turn to the proof of Theorem~\ref{T:critical_simple_symmetric}.
\newline {\sc Proof.} It follows from equation (\ref{P0(M>x)_main_equation}), Lemma~\ref{L:P_0(max>x)_simple_symmetric} and equality $1-f(1-s)=f'(c)s$, valid for $s\in[0,1]$ and some $c\in(1-s,1)$, that
$$
{\sf P}_0\left(M>x\right)\leq\left(\alpha f'(c)+(1-\alpha)\right){\sf P}_0\left(M>x\right)+{\sf P}_0\left(\max\left\{S(t),0\leq t\leq\tau_{0,0}\right\}>x\right).
$$
For $x\in\mathbb{Z}$ large enough, the latter inequality is possible only in the case when ${\sf P}_0\left(M>x\right)\to 0$, as $x\to\infty$.
Then according to the Taylor formula we have
\begin{eqnarray}\label{Taylor_critical}
& &1-f(1-{\sf P}_0\!\left(M>x\right))\\
&=&f'(1){\sf P}_0\!\left(M>x\right)-\frac{f''(1)}{2}\left({\sf P}_0\!\left(M>x\right)\right)^2+o\left(\left({\sf P}_0\!\left(M>x\right)\right)^2\right)\nonumber.
\end{eqnarray}
Hence, it follows from equation~(\ref{P0(M>x)_main_equation}) that
\begin{eqnarray}\label{P_0(M>x)^2(1+o(1))=}
& &\frac{\alpha\sigma^2}{2}\left({\sf P}_0\!\left(M>x\right)\right)^2(1+o(1))\\
&=&(1-\alpha){\sf P}_0\!\left(\max\left\{S(t),0\leq t\leq\tau_{0,0}\right\}>x\right)(1+o(1)),\nonumber
\end{eqnarray}
as $x\to\infty$.

Relations (\ref{P_0(max>x)_simple_symmetric}) and (\ref{P_0(M>x)^2(1+o(1))=}) imply the statement of Theorem~\ref{T:critical_simple_symmetric}. $\square$

Now let us turn to the proof of Theorem~\ref{T:subcritical_simple_symmetric}.
\newline {\sc Proof.} With the help of the same arguments as in the proof of Theorem~\ref{T:critical_simple_symmetric}, we conclude that ${\sf P}_0\left(M>x\right)\to 0$, as $x\to\infty$. However, in a subcritical case we write the Taylor formula in the form
\begin{equation}\label{Taylor_subcritical}
1-f(1-{\sf P}_0\!\left(M>x\right))=f'(1){\sf P}_0\!\left(M>x\right)+o\left({\sf P}_0\!\left(M>x\right)\right)
\end{equation}
and, reasoning in the same manner as in the proof of Theorem~\ref{T:critical_simple_symmetric}, we get
$$
\alpha\left(1-m\right){\sf P}_0\left(M>x\right)(1+o(1))=(1-\alpha){\sf P}_0\left(\max\left\{S(t),0\leq t\leq\tau_{0,0}\right\}>x\right)(1+o(1)),
$$
as $x\to\infty$. Whence the statement of Theorem~\ref{T:subcritical_simple_symmetric} follows. $\square$

Recall that $r$ is a probability of non-returning to point $0$ of the random walk $S$ starting from $0$. Let us prove Theorem~\ref{T:critical_simple_asymmetric}.
\newline {\sc Proof.} If $p<q$, then in view of Lemma~\ref{L:P_0(max>x)_simple_asymmetric} the following relations hold
\begin{equation}\label{P_0(S(t)<=x)_simple_asymmetric_sim_p<q}
{\sf P}_0\left(\max\left\{S(t),0\leq t\leq\tau_{0,0}\right\}\leq x,\tau_{0,0}<\infty\right)\to 2p,
\end{equation}
\begin{equation}\label{P_0(max>x)_simple_asymmetric_sim_p<q}
{\sf P}_0\left(\max\left\{S(t),0\leq t\leq\tau_{0,0}\right\}>x\right)\sim(q-p)\left(\frac{p}{q}\right)^{x+1},
\end{equation}
as $x\to\infty$. Moreover,
$$
{\sf P}_0\left(\max\left\{S(t),0\leq t\leq\tau_{0,0}\right\}\leq x,\tau_{0,0}<\infty\right)\to{\sf P}_0\left(\tau_{0,0}<\infty\right)=1-r,\quad x\to\infty.
$$
Consequently, $r=1-2p$ for $p<q$. By virtue of equation (\ref{P0(M>x)_main_equation}), formulae~(\ref{P_0(S(t)<=x)_simple_asymmetric_sim_p<q}), (\ref{P_0(max>x)_simple_asymmetric_sim_p<q}) and equality $1-f(1-s)=f'(c)s$, valid for $s\in[0,1]$ and some $c\in(1-s,1)$, we have
$$
{\sf P}_0\!\left(M\!>\!x\right)\!\leq\!\left(\alpha f'(c)+(1\!-\!\alpha)(1\!-\!r)\right){\sf P}_0\!\left(M>x\right)+{\sf P}_0\!\left(\max\left\{S(t),0\leq t\leq\tau_{0,0}\right\}\!>\!x\right).
$$
It follows that ${\sf P}_0(M>x)\to0$, as $x\to\infty$. Exploiting equation~(\ref{P0(M>x)_main_equation}), relations~(\ref{P_0(S(t)<=x)_simple_asymmetric_sim_p<q}), (\ref{P_0(max>x)_simple_asymmetric_sim_p<q}) and the Taylor formula in the form (\ref{Taylor_critical}) once again, we come to statement~(\ref{P0(M>x)_asymptotics_critical_simple_asymmetric_p<q}).

If $p>q$, then Lemma~\ref{L:P_0(max>x)_simple_asymmetric} implies that
\begin{equation}\label{P_0(S(t)<=x)_simple_asymmetric_sim_p>q}
{\sf P}_0\left(\max\left\{S(t),0\leq t\leq\tau_{0,0}\right\}\leq x,\tau_{0,0}<\infty\right)\to 2q
\end{equation}
and
\begin{equation}\label{P_0(max>x)_simple_asymmetric_sim_p>q}
{\sf P}_0\left(\max\left\{S(t),0\leq t\leq\tau_{0,0}\right\}>x\right)\to p-q,
\end{equation}
as $x\to\infty$. Then statement~(\ref{P0(M>x)_asymptotics_critical_simple_asymmetric_p>q}) follows from equation~(\ref{P0(M>x)_main_equation}) and reasoning on the existence and uniqueness of solution to this equation present in the proof of Lemma~\ref{L:main_equations}. Theorem~\ref{T:critical_simple_asymmetric} is proved completely. $\square$

It only remains to give the proof of Theorem~\ref{T:subcritical_simple_asymmetric}.
\newline {\sc Proof.} Employing the same arguments, as in the beginning of the proof of Theorem~\ref{T:critical_simple_asymmetric}, we come to conclusion that ${\sf P}_0\left(M>x\right)\to0$, as $x\to\infty$. Then applying relations (\ref{P0(M>x)_main_equation}), (\ref{Taylor_subcritical}), (\ref{P_0(S(t)<=x)_simple_asymmetric_sim_p<q}) and (\ref{P_0(max>x)_simple_asymmetric_sim_p<q}), we get formula (\ref{P0(M>x)_asymptotics_subcritical_simple_asymmetric_p<q}).

Statement~(\ref{P0(M>x)_asymptotics_subcritical_simple_asymmetric_p>q}) follows from relations~(\ref{P0(M>x)_main_equation}), (\ref{P_0(S(t)<=x)_simple_asymmetric_sim_p>q}), (\ref{P_0(max>x)_simple_asymmetric_sim_p>q}) and reasoning on the existence and uniqueness of solution to equation~(\ref{P0(M>x)_main_equation}), present in the proof of Lemma~\ref{L:main_equations}. Theorem~\ref{T:subcritical_simple_asymmetric} is proved completely. $\square$

In conclusion let us remark on the general case of an arbitrary finite number of catalysts in the critical CBRW on $\mathbb{Z}$. To investigate the asymptotic behavior of the solution to the system of equations (\ref{P0(M>x)_main_system}) we have to implement equivalent transformations of the system according to Cramer's rule (see, e.g., \cite{Kurosh_84}, Ch.~1, \S7), resulting in that the coefficient before ${\sf P}_{w_i}\left(M>x\right)$ for each $i=1,\ldots,N$ is equal to the determinant of the matrix $D-I$, where the matrix $D$ is specified in the definition of the critical regime and $I$ is the identity matrix. However, in the critical case $\det(D-I)=0$. Therefore, as in the case of a single catalyst, all the linear terms are reduced and there are quadratic terms only including $\left({\sf P}_{w_i}\left(M>x\right)\right)^2$. Other differences in the study of the solutions to equation (\ref{P0(M>x)_main_equation}) and system of equations (\ref{P0(M>x)_main_system}) are insignificant and for this reason we do not discuss them.

\vskip0.2cm The author expresses acknowledgements to Professors V.A.Vatutin, V.A.Topchij and S.G.Foss for useful discussions.

\end{document}